\documentclass[12pt]{article}

\usepackage{amsmath,amsthm,amsfonts,amssymb,bm }

\makeatletter \topmargin=0mm\headheight=0mm\headsep=0mm
\newcommand{\N}{{\mathbb N}}
\newcommand{\be}{\begin{eqnarray}}
\newcommand{\ee}{\end{eqnarray}}
\newcommand{\ben}{\begin{eqnarray*}}
\newcommand{\enn}{\end{eqnarray*}}

\newcommand{\sgn}{{\rm sgn}\,}

\newtheorem{theorem}{\textbf Theorem}[section]
\newtheorem{lemma}{\textbf Lemma}[section]

\def\endProof{{\hfill$\Box$}}

\newcounter{remark}
\newenvironment{remark}%
{\par \stepcounter{remark} {\it Remark
\arabic{section}.\arabic{remark}.}~}{\rm \endProof\par}

\def\R{\mathbb{R}}
\def\Da{\Lambda^\alpha}
\newcommand{\rf}[1]{{\rm (}\ref{#1}{\rm )}}
\allowdisplaybreaks

\begin{document}
\begin{titlepage}
\title{\bf On convergence of solutions of  fractal Burgers equation
toward rarefaction waves}
\author{Grzegorz Karch$^1$, Changxing Miao$^2$ and Xiaojing Xu$^2$\\
        \\
\small{      $^1$ Instytut Matematyczny, Uniwersytet Wroc{\l}awski,}\\
\small{ pl. Grunwaldzki 2/4, 50-384 Wroc{\l}aw, Poland}\\
\small{          (karch@math.uni.wroc.pl)}\\
\small{       $^{2}$ Institute of Applied Physics and Computational
Mathematics,}\\
\small{        P.O. Box 8009, Beijing 100088, P.R. China.}\\
\small{        (miao\_{}changxing@iapcm.ac.cn and
xjxu@email.jlu.edu.cn )}}

\date{}
\end{titlepage}
\maketitle

\begin{abstract}
In the paper, the large time  behavior of solutions
of the Cauchy problem for
 the one dimensional fractal Burgers equation
$u_t+(-\partial^2_x)^{\alpha/2} u+uu_x=0$ with $\alpha\in (1,2)$
 is studied.
It is shown that if the nondecreasing
initial datum approaches the  constant states $u_\pm$ ($u_-<u_+$) as
$x\to \pm\infty$, respectively, then the corresponding solution
converges toward the rarefaction wave, {\it i.e.}  the unique entropy
solution of the Riemann problem for the nonviscous Burgers equation.

\medskip

\noindent {\bf AMS Subject Classification 2000:}\quad 60J60, 35B40,
35K55,
35Q53.

\medskip

\noindent {\bf Key words:} fractal Burgers equation,
asymptotic behavior, rarefaction wave, Riemann problem, L\'evy process

\medskip

{\bf
\noindent
\date{\today}}

\end{abstract}

\section{Introduction}
\setcounter{equation}{0}

The goal of this work is to study asymptotic properties of solutions
to the Cauchy problem for the nonlocal conservation law
\begin{eqnarray}\label{e1.1}
&&u_t+\Da u +uu_x=0, \quad x\in\R,\; t>0,\\
\label{e1.2} &&u(0,x)=u_0(x)
\end{eqnarray}
where
$
\Da=(-\partial^2/\partial x^2)^{\alpha/2}
$
is the pseudodifferential operator defined {\it via} the Fourier
transform
\begin{equation}\label{Da}
\widehat {(\Da v)}(\xi)= |\xi|^\alpha \widehat v(\xi).
\end{equation}

Following \cite{bfw}, we will call \rf{e1.1} the fractal Burgers
equation.
Equations of this
type appear in the study of growing interfaces in the presence of
selfsimilar hopping
surface diffusion \cite{MW}. Moreover, in their recent papers,
Jourdain, M\'el\'eard, and Woyczynski
\cite{JMW, JMW1} gave probabilistic motivations to study equations
with
the anomalous diffusion, when Laplacian (the generator of the Wiener
process) is replaced by a more
 general pseudodifferential operator generating the L\'evy process.
In particular, the authors of \cite{JMW1}
studied problem \rf{e1.1}-\rf{e1.2}, where
the initial condition
$u_0$ is assumed to be a nonconstant function with bounded variation
on $\R$. In
other words, a.e. on $\R$,
\begin{equation}\label{u0}
u_0(x)=c+\int_{-\infty}^xm(dy)=c+H*m(x)
\end{equation}
with $c \in \R$, $m$ being a finite signed measure on $\R$, and $H(y)$
denoting
the unit step function $1\!\!{\rm I}_{\{y\geq 0\}}$.
Observe that the gradient $v(x,t)=u_x(x,t)$ satisfies
\begin{equation}\label{FP}
v_t+\Da v +(vH*v)_x=0, \quad v(\cdot, 0)=m.
\end{equation}
If m is a probability measure on $\R$, the equation \rf{FP} is a
nonlinear Fokker-Planck
equation. In the case of an arbitrary finite signed measure,
 the authors of  \cite{JMW1}
 associate \rf{FP} with a suitable nonlinear martingale problem.
 Next, they study the convergence of systems of particles with jumps
as the number of particles tends to $+\infty$.
 As a consequence, the weighted empirical cumulative distribution
functions
 of the particles converge to the solution of the martingale problem
connected to \rf{FP}. This phenomena is called
 {\it the propagation of chaos} for problem \rf{e1.1}-\rf{e1.2} and we
refer the reader to \cite{JMW1} for
 more details and additional references.

Motivated by the results from \cite{JMW1}, we study problem
\rf{e1.1}-\rf{e1.2}
 under   the crucial assumption
$\alpha\in (1,2)$ and with the
 initial condition of the form \rf{u0}. In our main result, we assume
that
  $u_0$
  is a
 function satisfying
\begin{equation}\label{ini as}
u_0-u_-\in L^1((-\infty,0)) \quad\mbox{and}\quad
 u_0-u_+\in L^1((0,+\infty))
\quad \mbox{with}\quad u_-<u_+,
\end{equation}
where $u_-=c$ and $u_+=c+\int_{\R}m(dx)$.

It is well known (cf. \cite{ilo,man,han} and Lemma \ref{le3.6}, below)
that the asymptotic profile as $t\to\infty$  of
 solutions to the viscous Burgers equation
\begin{equation}\label{burgers}
u_t-u_{xx}+uu_x=0
\end{equation}
 ({\it i.e.} equation (\ref{e1.1}) with $\alpha=2$) supplemented with
an
 initial datum satisfying \rf{ini as} is given
by the so-called rarefaction wave. This is the continuous function
\begin{equation}\label{rarefaction}
w^R(x,t)=W^R(x/t)=\left\{
\begin{aligned}
&u_-\,,\quad &&x/t\leq u_-,\\
&x/t\,,\quad &&u_-\leq x/t\leq u_+, \\
& u_+\,,\quad &&x/t\geq u_+,
\end{aligned}
\right.
\end{equation}
which is the unique entropy solution of
following Riemann problem
\begin{eqnarray}
&&w^R_t+w^R w^R_x=0,\label{eq-rar}\\
&&w^R(x,0)=w_0^R(x)=\left\{
\begin{aligned}
u_-\,,\quad x<0,\\
u_+\,,\quad x>0.
\end{aligned}\right.\label{ini-rar}
\end{eqnarray}
Below, we use the solution of the Burgers equation \rf{burgers} with
the initial datum \rf{ini-rar} as the smooth approximation of the
rarefaction wave \rf{rarefaction}.

The authors of this work were inspired by the fundamental paper of
 Il'in
and Oleinik \cite{ilo} who showed the  convergence  toward
rarefaction waves of solutions to the nonlinear equation
$u_t-u_{xx}+f(u)_x=0$ under strict convexity assumption
imposed on $f$.
That idea was   next extended in several different directions and we
refer the reader
{\it e.g.} to \cite{han,lmn,man,man1,nak,nin}
for an overview of know results and additional references.

In this work, we contribute to the existing theory by developing tools
which allows to obtain analogous results for equations with a nonlocal
and anomalous diffusion.
Basic properties of solutions (namely, their existence and the
regularity) of
quasilinear evolution equations with 
$(-\Delta)^{\alpha/2}$,
$\alpha\in (1,2)$, (or, more generally, with the L\'evy diffusion)
were shown in \cite{dgv,dri}.
On the other hand, one may expect singularities in finite time 
of solutions to \rf{e1.1} with $\alpha\in (0,1)$, see \cite{ADV} for more details. 
If  $u_0\in L^1(\R)$, Biler, Karch, and
Woyczynski \cite{bkw} proved that the large time asymptotics of
solutions to
\rf{e1.1}-\rf{e1.2} is described by the self-similar fundamental
solution
of equation $v_t+\Da v=0$. Analogous asymptotic properties of
solutions to multidimensional
generalizations of problem \rf{e1.1}-\rf{e1.2} with $u_0\in L^1(\R^N)$
are studied in \cite{bkw1,bkw2}.

Finally, we would like to report the recent progress in the
understanding of properties of solutions
 of the quasi-geostrophic equation
with an anomalous diffusion, {\it cf.} \cite{cc,ju} and the references
therein.
We are convinced that our techniques can be applied to that equation,
as well.

\bigskip

The purpose of the present paper is to
prove the convergence of  solutions to the Cauchy problem for the
fractal Burgers equation
\rf{e1.1}--\rf{e1.2} toward rarefaction waves.
We state our main result  in the following theorem.

\begin{theorem}\label{th main}
Let $\alpha\in (1,2)$.
Assume that $w^R=w^R(x,t)$ is the rarefaction wave \rf{rarefaction}.
Denote by $u=u(x,t)$ the unique solution of problem
\rf{e1.1}-\rf{e1.2}
corresponding to the initial datum $u_0$ of the form \rf{u0} and
satisfying \rf{ini as}.
For every $p\in ((3-\alpha)/(\alpha-1), \infty]$ there exists
$C>0$  independent of $t$ such that
\begin{align*}
\|u(t)-w^R(t)\|_{p}\leq
Ct^{-[\alpha-1-(3-\alpha)/p]/2}\log (2+t)
\end{align*}
for all $t>0$. 
\end{theorem}

\begin{remark}
Our result and its proof hold true also for $\alpha=2$
(observe that $(3-\alpha)/(\alpha-1)\to 1$ as $\alpha\to 2$). However,
we pass over this case for simplicity of the exposition and because
the large time asymptotics of solutions to the Burgers equation
\rf{burgers} is well-known, see Lemma \ref{le3.6}, below.
\end{remark}

In the next section, we gather several preliminary properties of the
operator $\Da$
and of  solutions to problem \rf{e1.1}-\rf{e1.2}. Theorem \ref{th
main} is shown in Section 3.
In Section 4, we discuss possible generalizations of our main result.

{\bf Notation.}
For
$1\leq p\leq \infty$, the $L^p$-norm   of a Lebesgue
measurable, real-valued function $v$ defined on $\R$
is denoted by $\|v\|_p$.
For a finite signed measure $m$ on $\R$, we put $\|m\|=|m|(\R)$, where
$|m|$ is the total variation of $m$.
The Fourier transform  of $v$ is
$\widehat v(\xi)\equiv (2\pi)^{-1/2}\int_{\R}
e^{-ix\xi} v(x)\;dx$.
Given a function $v=v(x)$, we are going to use the decomposition
$
v=v^-+v^+,
$
where as usual  $v^-=\max\{0,\,-v\}$ and $v^+=\max\{0,\,v\}.$
The constants (always independent of
of $t$) will be
denoted by the
same letter $C$, even if they may vary from line to line.
Occasionally, we write,  e.g.,  $C=C(\alpha,\ell)$ when we want to
emphasize
the dependence of $C$ on parameters $\alpha$ and $\ell$.

\section{Preliminary results}
\setcounter{equation}{0}

We begin  by recalling that the basic questions on
the existence and the uniqueness of solutions of  problem
\rf{e1.1}-\rf{e1.2} were answered  in the papers  \cite{dgv,dri}.

\begin{theorem}\label{th2.1} {\rm (\cite[Thm. 1.1]{dgv}, \cite[Thm.
7]{dri})}
Let $\alpha\in (1,\,2)$ and $u_0\in L^\infty(\R)$.
There exists the unique solution $u=u(x,t)$ to problem
\rf{e1.1}-\rf{e1.2} in
the following sense: for all $T>0$,
\begin{align*}
&u\in C_b((0,\,T)\times\mathbb{R})\;and,\;
for\;all\;a\in(0,\,T),\;u\in C^\infty_b((a,\,T)\times\mathbb{R}),
\\
&u\;satisfies\;(\ref{e1.1})\;on \;(0,\,T)\times\mathbb{R},
\\
&u(t,\cdot)\rightarrow u_0\;in
\;L^\infty(\mathbb{R})\;weak-\ast,\;as \; t\rightarrow 0.
\end{align*}
Moreover, the following inequality holds true
\begin{equation}
\|u(t)\|_\infty\leq \|u_0\|_\infty \quad \mbox{for all}\quad t>0.
\end{equation}
\end{theorem}

The main goal of the section is to complete this result by additional
properties
of $u_x$  if the initial conditions are of the form \rf{u0}.

\begin{theorem}\label{th ux}
Let $\alpha\in (1,\,2)$.
Assume that the initial datum $u_0$ can be written in the form \rf{u0} for a
constant $c\in \R$
and a signed finite measure $m$ on $\R$.
Then  the solution $u=u(x,t)$ of problem
 \rf{e1.1}-\rf{e1.2} satisfies $u_x\in C((0,T];L^p(\R))$  for each $1\le p\le\infty$ and every
$T>0$.

Consider  $u$ and $\widetilde u$ two such  solutions  with initial conditions $u_0$ and $\widetilde
 u_0$, respectively. Suppose  that
 $\widetilde u_x(x,t)$ is nonnegative a.e.  and
 $u_0-\widetilde u_0\in L^1(\R)$.
Then 
\begin{equation}
\|u(t)-\widetilde u(t)\|_1\leq \|u_0-\widetilde u_0\|_1\label{L1}
\end{equation}
for all $t>0$.
\end{theorem}

\begin{theorem}\label{th ux 2}
Under the assumption of Theorem \ref{th ux},  
if the measure $m$ in the initial datum \rf{u0} is nonnegative,  we have

i) $u_x(x,t)\geq 0$ for all $x\in\R$ and $t>0$,

ii) for every $p\in [1,\infty]$ there exists $C=C(p)>0$ such that
\begin{equation}\label{ux dec}
\|u_x(t)\|_p\leq t^{-1+1/p}\|m\|^{1/p}.
\end{equation}
\end{theorem}

In the proofs of Theorems \ref{th ux} and \ref{th ux 2}
as well as in our study of the large time asymptotics,
we shall  require several properties of the
operator $\Da$ and of the semigroup of linear operators generated by
it.
First of all, note that the operator defined by \rf{Da} has
the
integral representation for every $\alpha\in (1,2)$
(cf. eg.  \cite[Thm. 1]{dri})
\begin{equation}\label{rep}
\Da w(x)=-C(\alpha) \int_{\R}\frac{w(x+z)-w(x)- w_x(x)
z}{|z|^{1+\alpha}}\;dz.
\end{equation}
This formula allows us to apply $\Da$ to functions which are bounded
and sufficiently smooth, however, not necessary decaying at infinity.

\bigskip

\begin{lemma}
Let $1<\alpha<2$.
For every $p\in[1,\,\infty]$ there exists $C=C(p,\alpha)>0$ such that
\begin{equation}\label{interpol}
\|\Da w\|_p\leq C\|w_x\|^{2-\alpha}_p\|w_{xx}\|^{\alpha-1}_p
\end{equation}
all functions $w$ satisfying  $w_x,w_{xx}\in L^p(\R)$.
\end{lemma}

\begin{proof}
We can easily deduce the interpolation inequality \rf{interpol} from
\rf{rep}.
Indeed, it follows from  the Taylor formula that for any fixed $R>0$
we have
\begin{eqnarray*}
\|\Da w\|_p&\leq& C\|w_{xx}\|_p\int_{|z|\leq R} |z|^{1-\alpha}\;dz
+C\|w_x\|_p\int_{|z|>R} |z|^{-\alpha} \;dz\\
&\leq& C\left( R^{2-\alpha}\|w_{xx}\|_p +R^{1-\alpha}\|w_x\|_p\right).
\end{eqnarray*}
Choosing $R=\|w_x\|_p/\|w_{xx}\|_p$ we complete the proof of
inequality \rf{interpol}.
\end{proof}

Now, we prove the Nash inequality for the operator $\Da$.

\begin{lemma}
Let $0<\alpha$.
There exists a constant $C_N>0$ such that
\begin{equation}\label{Nash}
\|w\|_2^{2(1+\alpha)} \leq C_N
\|\Lambda^{\alpha/2}w\|_2^2\|w\|_1^{2\alpha}
\end{equation}
for all functions $w$ satisfying  $w \in L^1(\R)$ and
$\Lambda^{\alpha/2}w\in L^2(\R)$.
\end{lemma}

\begin{proof}
For every $R>0$, we decompose the $L^2$-norm of the Fourier transform
of $w$ as follows
\begin{align*}
\|w\|_2^2&=C\int_\R |\widehat w(\xi)|^2\;d\xi\\
&\leq C\|\widehat w\|_\infty^2\int_{|\xi|\leq R}\,d\xi
+CR^{-\alpha}\int_{|\xi|>R} |\xi|^\alpha|\widehat w(\xi)|^2\,d\xi\\
&\leq CR \|w\|_1^2+CR^{-\alpha}\|\Lambda^{\alpha/2}w\|_2^2.
\end{align*}
For
$R=\left(\|\Lambda^{\alpha/2}w\|_2^2/\|w\|_1^2\right)^{1/(1+\alpha)}$
we obtain \rf{Nash}.
\end{proof}

\begin{lemma}\label{le3.1}
Let $0\leq \alpha\leq 2$.
For every
$p>1$,  we have
\begin{equation}\label{hyper}
\int_\R (\Da w)|w|^{p-2} w \,dx\geq
\frac {4(p-1)}{p^2}\int_\R\left(\Lambda^{\frac\alpha
2}|w|^{\frac p2}\right)^2\;dx
\end{equation}
for all $w\in L^p(\R)$ such that $\Da w\in L^p(\R)$. If $\Da w\in
L^1(\R)$, we obtain
\begin{equation}\label{sgn est}
\int_\R (\Da w)\sgn w\,dx\geq 0,
\end{equation}
and, if $w,\Da w\in L^2(\R)$, it follows
\begin{equation}\label{+ est}
\int_\R (\Da w) w^+\,dx\geq 0 \quad\mbox{and}\quad \int_\R (\Da w)
w^-\,dx\geq 0,
\end{equation}
where $w^+=\max\{0,w\}$  and $w^-=\max\{0,-w\}$.
\end{lemma}

Inequality \rf{hyper} is well-known in the theory of sub-Markovian
operators and its statement
 and the proof is given {\it e.g.} in
\cite[Theorem 2.1 combined with the Beurling-Deny condition
(1.7)]{LS},
see also \cite{cc,ju}.
Observe that if $\alpha=2$, integrating by parts we obtain \rf{hyper}
with the equality.
Inequality \rf{sgn est} (called the Kato inequlity)
is used in \cite{dri} to construct entropy solutions of \rf{e1.1}
and it can be easily deduced from \cite[Lemma 1]{dri} by an
approximation argument (see also \cite[Inequality (3.5)]{bfw}).
The proof of \rf{+ est} can be found, for example, in
\cite[Proposition 1.6]{LS}.

We also recall that,
by Duhamel's principle, the solution to problem \rf{e1.1}-\rf{e1.2}
can be written in the
equivalent integral form
\be
u(t)=S_{\alpha}(t)u_{0}-\int^t_0S_{\alpha}(t-\tau)u(\tau)u_x(\tau)\,d\
tau,
\label{duhamel}
\ee
where
\begin{equation}\label{e3.10}
S_{\alpha}(t)u_{0}=p_\alpha(t)*u_{0}(x).\end{equation}
Here, the fundamental solution
$p_\alpha(x,t)$
of the linear equation
$ \partial_t v+\Da v=0$
can be computed {\it via} the Fourier
transform
$ \widehat p_\alpha(\xi,t)=e^{-t|\xi|^\alpha}$.
Hence,
$
p_\alpha(x,t)=t^{-1/\alpha}P_\alpha(xt^{-1/\alpha}),
$
where $P_\alpha$ is the inverse Fourier transform of
$e^{-|\xi|^\alpha}$.
It is well known that for every $\alpha\in (0,2]$ the function
$P_\alpha$
has the property $\int_\R P_\alpha(x)\;dx=1$ and
is  smooth, nonnegative,  and
satisfies 
\begin{equation}
\label{EP}
0\le P_\alpha(x)\le C(1+|x|)^{-(\alpha+1)} \qquad \hbox{and}\qquad
|\partial_x P_\alpha(x)|\le C(1+|x|)^{-(\alpha+2)}
\end{equation}
for a constant $C$ and all $x\in\R$.
Using these properties of the convolution operator $S_{\alpha}(t)$
defined by (\ref{e3.10})
we obtain the estimates
\be\label{e3.11}
\|S_\alpha(t)v\|_p&\le& Ct^{-(1-1/p)/\alpha)}\|v\|_{1},\\
\label{e3.12}
\| (S_\alpha(t)v)_x\|_p
&\le& Ct^{-(1-1/p)/\alpha-1/\alpha}\|v\|_{1}
\ee
for every $p\in [1,\infty]$ and all $t>0$.
Moreover, we can replace $v$ in \rf{e3.11} and in \rf{e3.12} by any
signed measure $m$.
In that case, $\|v\|_1$ should be replaced by $\|m\|$.

\bigskip
\noindent
{\it Proof of Theorem \ref{th ux}.}
It follows from the integral equation \rf{duhamel} that $u_x$ is the
solution of
\begin{equation}
u_x(t)=S_\alpha(t)m-\int_0^t \partial_x S_\alpha(t-\tau)
V(\tau)u_x(\tau)\,d\tau, \label{duhamel dx}
\end{equation}
where $V(x,t)=u(x,t)$ is treated as given, smooth, and bounded.
Now the standard argument involving the Banach fixed point theorem
allows us to show that
the solution of the ``linear'' equation \rf{duhamel dx} has the
solution
in $C((0,T];L^p(\R))$ for each $p\in [1,\infty]$ and every $T>0$.
Here, we should
use the following estimate
of the operator ${\cal T} (u)$ defined by the right-hand side of
\rf{duhamel dx}
\begin{align*}
\|{\cal T}(u)(t)\|_p&\leq
\|S_\alpha(t)m\|_p+\int_0^t \|\partial_x S_\alpha(t-\tau)
V(\tau)u_x(\tau)\|_p\,d\tau\\
&\leq Ct^{-(1-1/p)/\alpha}\|m\|
+ C \sup_{\tau\in [0,T]}\|V(\tau)\|_\infty\int_0^t
(t-\tau)^{-1/\alpha} \|u_x(\tau)\|_p\,d\tau,
\end{align*}
being the immediate consequence of \rf{e3.11} and \rf{e3.12}. Let us
skip other details
of this well-known argument (cf. \cite{MYZ}).

Now, we prove inequality \rf{L1}.
A direct calculation shows that the function $v(x,t)=u(x,t)-\widetilde
u(x,t)$ satisfies
\begin{equation}\label{eq v}
v_t+\Da v +\frac12 (v^2+2v \widetilde u))_x=0.
\end{equation}
First, we multiply equation \rf{eq v} by $\sgn v=v|v|^{-1}$:
\begin{equation*}
\frac d{dt}\int_\R|v|\,dx+
\int_\R (\Da v) \sgn v\,dx+\frac 12\int_\R[v^2+2v\widetilde u)]_x \sgn 
v
\,dx=0.
\end{equation*}
The second term is nonnegative by \rf{sgn est}.
To show the same property for the third term, we replace the
sgn function  by  smooth and nondecreasing
$\varphi=\varphi(x)$.
In this case, we obtain
$$
\int_{\R} [v^2+2v\widetilde u]_x \varphi(v) dx
=-\int_\R (v^2+2v\widetilde u)\varphi'(v)v_x \;dx
=-   \int_\R \Psi(v)_x\; dx+
\int_\R \widetilde u_x \Phi(v)\; dx,
$$
where $\Psi(s)=\int_0^sz^2\varphi'(z)\,dz$ and
$\Phi(s)=\int_0^s2z\varphi'(z)\,dz$.
Obviously, the first term on the right hand side is equal to zero and
the second one is nonnegative
because  $\widetilde u_x\geq 0$ and
$\Phi(s)\geq 0$ for all $s\in \R$.
Now, the standard approximation argument gives
$\int_\R[v^2+2v\widetilde u]_x \sgn v \,dx\geq 0$. Hence
$
\|v(t)\|_1=\|u(t)-\widetilde u(t)\|_1\leq \|u_0-\widetilde u_0\|_1
=\|v_0\|_1
$
for all $t>0$.
\endProof

\bigskip
\noindent
{\it Proof of Theorem \ref{th ux 2}.}
To show part (i) of Theorem \ref{th ux 2}, we deal
first with the smooth initial datum $u_0$ satisfying $u_{0,x}(x)\geq 
0$ and $u_{0,x}\in L^p(\R)$ for every $p\in [1,\infty]$.
In this case, differentiating equation (\ref{e1.1}) with respect to 
$x$ 
we 
have
\begin{align}\label{e3.1}
(u_x)_t+\Da u_x+(uu_x)_x=0.
\end{align}
Note the well known property
\begin{align*}
\int_\R v_tv^-\,dx=\int_{v\leq 0} v^-_tv^-\,dx=\frac12\frac
d{dt}\int_{v\leq 0}(v^-)^2\,dx.
\end{align*}
Hence,  multiplying (\ref{e3.1}) by $u_x^-$, integrating the resulting
equation
over $\R$, and integrating by parts on the right hand side, we obtain
\begin{align*}
\frac12\frac
d{dt}\int_{u_x\leq 0}(u_x^-)^2\;dx
+
\int_{\R}(\Da u_x) u_x^-\;dx
&=-\int_{u_x\leq 0}(uu_x^-)_xu_x^-\;dx
\\
&=-\frac12\int_{u_x\leq 0}(u_x^-)^3\,dx.
\end{align*}
Since $\int_{\R}(\Da u_x) u_x^-\,dx\geq 0$ by \rf{+ est} and
$\int_{u_x\leq 0}(u_x^-(x,0))^2\,dx=0$
by the assumption imposed on $u_0$,  the
Gronwall inequality implies
$
\int_{u_x\leq 0}(u_x^-(x,t))^2\,dx= 0
$
for all $t\geq 0$.
Consequently,
$
u_x^-(x,t)\equiv 0
$
and  the proof of (i) for regular initial conditions is finished.

Now, the proof of (i) for the solution $u=u(x,t)$ corresponding to
 the initial datum $u_0$ of the form \rf{u0}
with the nonnegative finite measure $m$ can be completed by the following approximation argument.
We consider the sequence of regular initial conditions $u_0^n$ as in the first part of this proof.
Moreover, we assume that $u_{0,x}^n$ converges weakly to $m$ and $\|u_0^n-u_0\|_1\to 0$  as $n\to\infty$.
Inequality \rf{L1}  allows us to prove that the corresponding 
solutions $u^n(\cdot, t)$ satisfy $\|u^n(\cdot,t)-u(\cdot,t)\|_1\to 0$ as $n\to\infty$
 for any $t>0$. 
Hence, there is a subsequence $n_k\to \infty$ such that $u^{n_k}(x,t)\to u(x,t)$ a.e.
Since each $u^n(x,t)$ is nondecreasing as  function of $x$, the same conclusion holds true for $u(x,t)$.

In order to show inequality \rf{ux dec}, we first observe that integrating
equation \rf{duhamel dx}
over $\R$ 
and using the equalities
$$
\int_\R S_\alpha(t)m\;dt=\int_\R m(dx)\quad\mbox{and} \quad 
\int_\R \partial_x S_\alpha(t-\tau)(u(\tau)u_x(\tau))\;dx=0
$$
we obtain the identity $\int_\R u_x(x,t)\,dx =\int_\R m(dx)$
which
for nonnegative
$u_x$ means
\begin{equation}\label{ux L1}
\|u_x(t)\|_1=\|m\| 	\quad \mbox{for all} \quad t>0.
\end{equation}

Now, for fixed $p\in (1,\infty)$ and $u_x\geq 0$, we multiply
\rf{e3.1} by $u_x^{p-1}$ and integrate the resulting equation over
$\R$.
After some manipulations involving integrations by parts on the right
hand side
we arrive at
\begin{align}\label{ux est}
\frac1p\frac d{dt}\|u_x\|^p_p+\int_\R u_x^{p-1} \Da u_x\,dx&=-\int_\R
(uu_x)_xu_x^{p-1}\,dx \nonumber \\
&=-\frac{p-1}p\int_\R u_x^{p+1}\,dx.
\end{align}
Recall now that $\int_\R u_x^{p-1} \Da u_x\,dx\geq 0$ by inequality
\rf{hyper}.
Moreover, it follows from the H\"older inequality combined with \rf{ux
L1}
that
$$
\|u_x(t)\|_p^{p^2/(p-1)} \leq \|u_x(t)\|_{p+1}^{p+1}\|m\|^{1/(p-1)}.
$$
Applying those two inequalities to \rf{ux est} (note $u_x\geq 0$)
we obtain the following differential inequality for
$\|u_x(t)\|_p^p$
$$
\frac{d}{dt} \|u_x(t)\|_p^p \leq -(p-1) \|m\|^{-1/(p-1)}
\left(\|u_x(t)\|_p^p\right)^{p/(p-1)}.
$$
Integrating it we complete the proof of \rf{ux dec} for any $p\in
(1,\infty)$.

The case of $p=\infty$ is obtained immediately passing to the limit
$p\to \infty$ in inequality~\rf{ux dec}.
\endProof

\bigskip

We conclude this section by recalling some results on smooth
approximations of
rarefaction waves, namely, the solutions
 of the following Cauchy problem
\begin{eqnarray}
&&w_t-w_{xx}+w w_x=0,\label{eq-app}\\
&&w(x,0)=w_0(x)=
\left\{
\begin{aligned}
u_-\,,\quad x<0,\\
u_+\,,\quad x>0.
\end{aligned}\right.
.
\label{ini-app}
\end{eqnarray}

\begin{lemma}\label{le3.6}
Let $u_-<u_+$.
Problem \rf{eq-app}-\rf{ini-app}  has the unique, smooth,
global-in-time  solution $w(x,t)$ satisfying

i) $u_-<w(t,x)<u_+$ and  $w_x(t,x)>0$ for all
$(x,t)\in\R\times(0,\infty)$;

ii)  for every $p\in[1,\;\infty]$, there exists a constant
$C=C(p,u_-,u_+)>0$ such that
\begin{equation*}
\|w_x(t)\|_{p}\leq C t^{-1+1/p},
\quad
\|w_{xx}(t)\|_{p}\leq C t^{-3/2+1/(2p)}
\end{equation*}
and
\begin{equation*}
\|w(t)-w^R(t)\|_{p}\leq C t^{-(1-1/p)/2},
\end{equation*}
for all $t>0$, where $w^R(x,t)$ is the rarefaction wave
\rf{rarefaction}.
\end{lemma}

All results stated in Lemma \ref{le3.6} are deduced from the explicit
formula for
solutions to \rf{eq-app}-\rf{ini-app} and detailed calculations can be
found in
\cite{han} with some additional improvements contained in
\cite[Section 3]{kat}.

\section{Convergence toward rarefaction waves}
\setcounter{equation}{0}

For simplicity of the exposition, we divide the proof of Theorem
\ref{th main} into a sequence of Lemmata.

\begin{lemma}\label{lem 1}
Let $\alpha\in  (1,2)$.
 Assume that $u$ and $\widetilde u$ are two solutions of problem
 \rf{e1.1}-\rf{e1.2} with initial conditions $u_0$ and $\widetilde
 u_0$,  the both of the from \rf{u0} with finite signed
 measures $m$ and $\widetilde m$, respectively.
 Suppose,
 moreover, that
 the measure $\widetilde m$ of $\widetilde u_0$ is nonnegative and
 $u_0-\widetilde u_0\in L^1(\R)$. Then, for every $p\in
 [1,\infty]$ there exists a constant $C=C(p)>0$ such that
\begin{equation}\label{u tilde}
\|u(t)-\widetilde u(t)\|_p\leq Ct^{-(1-1/p)/\alpha}\|u_0-\widetilde
u_0\|_1
\end{equation}
for all $t>0$.
\end{lemma}

\begin{proof}
In our reasoning, we denote $v(x,t)=u(x,t)-\tilde u(x,t)$ which satisfies equation \rf{eq v}. 
It follows from Theorem \ref{th ux}, inequality \rf{L1}, that 
$\|v(t)\|_1\leq \|v_0\|_1$.

Now,
we multiply  equation \rf{eq v} by $|v|^{p-2}v$ with $p>1$
\begin{equation}\label{3.3a}
\frac 1p\frac d{dt}\int_\R|v|^p\,dx+
\int_\R (\Da v) (|v|^{p-2}v)\,dx+\frac 12\int_\R[v^2+2v\widetilde
u]_x|v|^{p-2}v\, dx
=0.
\end{equation}
The third term on the left hand side  of  \rf{3.3a}
is  nonnegative  by  the following calculations
\begin{align}\label{v pos}
\int_{\R} [v^2+2v\widetilde u]_x|v|^{p-2}v \,dx&=
\int_\R 2 v_x|v|^p\,dx+
\int_\R 2\widetilde u v_x|v|^{p-2}v\,dx+
\int_\R 2\widetilde u_x|v|^p \,dx
\\
&=2\left(1-\frac1p\right) \int_\R \widetilde u_x|v|^p\; dx \geq
0,\nonumber
\end{align}
because $\int_\R v_x|v|^p\,dx=0$ and $\widetilde u_x\geq 0$.
Hence, using inequality \rf{hyper}, we obtain from \rf{3.3a}
\begin{equation}\label{est p}
\frac d{dt}\int_\R|v|^p\,dx+4 \left(1-\frac1p\right)\int_\R
(\Lambda^{\alpha/2}
|v|^{p/2})^2\,dx\leq 0.
\end{equation}

From now on, we proceed by induction.
Applying the Nash inequality \rf{Nash} combined with \rf{L1}, we
deduce from \rf{est p}
with $p=2$ the following differential inequality
$$
\frac d{dt} \|v(t)\|_2^2 +
2C_N^{-1}\|v_0\|_1^{-2\alpha}\|v(t)\|_2^{2(1+\alpha)}\leq 0,
$$
which, after integration, leads to
\begin{equation}\label{p=2}
\|v(t)\|_2\leq C_1 \|v_0\|_1t^{-1/(2\alpha)}
\quad\mbox{with} \quad
C_1=( C_N/2\alpha)^{1/(2\alpha)}.
\end{equation}
This is estimate
\rf{u tilde} with $p=2$.

Suppose now that
\begin{equation}
\|v(t)\|_{2^n} \leq C_n t^{-(1-2^{-n})/\alpha}\|v_0\|_1 \label{ind as}
\quad\mbox{for all} \quad t>0.
\end{equation}
We consider  \rf{est p} with $p=2^{n+1}$, where the second term is
estimated, first, by the
Nash inequality \rf{Nash} with $w=|v|^{2^n}$, next, by the inductive
hypothesis \rf{ind as}.	
This two-step estimate leads to the differential inequality
$$
\frac d{dt} \|v(t)\|_{2^{n+1}}^{2^{n+1}} +
4(1-2^{-n-1})C_N^{-1} (C_n\|v_0\|_1)^{-\alpha2^{n+1}}
t^{2^{n+1}-2}\left(\|v(t)\|_{2^{n+1}}^{2^{n+1}}\right)^{1+\alpha}
\leq 0.
$$
Integrating it we obtain
\begin{equation}\label{ind as n+1}
\|v(t)\|_{2^{n+1}} \leq C_{n+1}
t^{-(1-2^{-n-1})/\alpha} \|v_0\|_1 \quad\mbox{for all} \quad
t>0,
\end{equation}
with
$$
C_{n+1}=C_n \left((C_N/(2\alpha))^{1/\alpha}\right)^{2^{-n-1}}
\left(
2^{n2^{-n-1}}
\right)^{1/\alpha}.
$$
This is  inequality \rf{u tilde} for
any $p=2^{n+1}$ with $n\in\N$.

 We leave to the reader the proof that $\lim\sup_{n\to \infty}
C_n<\infty$. Hence, passing
 to the limit $n\to\infty$ in \rf{ind as n+1} we obtain inequality
\rf{u tilde} for $p=\infty$.

The H\"older inequality
$$
\|v\|_p\leq \|v\|_{2^n}^{2^{n+1}/p -1}\|v\|_{2^{n+1}}^{2-2^{n+1}/p}
$$
completes the proof for every
$p\in (2^n, 2^{n+1})$.
\end{proof}

\begin{lemma}\label{lem 2}
Let $\alpha\in  (1,2)$.
Assume that $w=w(x,t)$ is the smooth approximation of the rarefaction
wave,
namely, the solution of problem \rf{eq-app}-\rf{ini-app}.
Then for each $t_0>0$ we have
$$
\int_{t_0}^\infty \| w_{xx}(t)\|_p\;dt <\infty
\quad\mbox{for every}
\quad
p \in (1, \infty]
$$
and
$$
\int_{t_0}^t \| \Da w(t)\|_p\;dt \leq C\log (2+t)
\quad\mbox{for }
\quad
p =\frac{3-\alpha}{\alpha-1},
$$
all $t\geq t_0$ and $C>0$ independent of $t$.
\end{lemma}

\begin{proof}
It follows from the decay estimates recalled in Lemma \ref{le3.6} that
$$
\int_{t_0}^\infty \| w_{xx}(t)\|_p\;dt \leq C
\int_{t_0}^\infty t^{-3/2+1/(2p)}\;dt
<\infty
\quad\mbox{for every}
\quad
p \in (1, \infty].
$$
By the interpolation inequality
\rf{interpol} and  Lemma \ref{le3.6}, we
obtain
\begin{eqnarray*}
\|\Da w(t)\|_p&\leq&
C(1+t)^{(-1+1/p)(2-\alpha)}(1+t)^{(-3/2+1/(2p))(\alpha -1)}\\
&=&
C(1+t)^{-(1+\alpha)/2+(3-\alpha)/(2p)}.
\end{eqnarray*}
Hence, the rate of decay on the right hand side is equal to $-1$
for 
$p=(3-\alpha)/(\alpha-1)$.
\end{proof}

\begin{lemma}\label{lem 3}
Let $\alpha\in  (1,2)$.
Assume that $u=u(x,t)$ is the solution of  \rf{e1.1}-\rf{e1.2}
and
$w=w(x,t)$ -- of  \rf{eq-app}-\rf{ini-app}.
Suppose that $u_0-w_0\in L^p(\R)$
for   $p=(3-\alpha)/(\alpha-1)$.
Then
\begin{equation*}
 \|u(t)-w(t)\|_p\leq C\log (2+t).
\end{equation*}
\end{lemma}

\begin{proof}
Denoting $v=u-w$, we see that this new function satisfies
\begin{equation*}
v_t+\Da v+\frac12[v^2+2vw]_x=-\Da w+w_{xx}.
\end{equation*}
We multiply this equation by $|v|^{p-2}v$ and we integrate over
$\R$ to obtain
\begin{equation}\label{e4.2}
\begin{aligned}
&\frac 1p\frac d{dt}\int|v|^pdx+
\int (\Da v) (|v|^{p-2}v)+\frac 12\int[v^2+2vw]_x|v|^{p-2}v dx
\\
=&\int (-\Da w+w_{xx})(|v|^{p-2}v)dx.
\end{aligned}
\end{equation}
It follows from Lemma \ref{le3.1} that
$
\int_{\R} (\Da v) (|v|^{p-2}v)dx\geq 0.
$
The third term on the left hand side  of (\ref{e4.2})
is  nonnegative  by  Lemma \ref{le3.6}
and the same argument as the one used in the proof of Lemma \ref{lem
1}, cf. identity \rf{v pos}.
Moreover, using the H$\ddot{\rm o}$lder inequality, we have
\begin{align*}
\left|\int_{\R} (-\Da w+w_{xx})(|v|^{p-2}v)\; dx\right|
\leq  (\|\Da w\|_p+\|w_{xx}\|_p)\|v\|_p^{p-1}.
\end{align*}

Consequently, \rf{e4.2} implies the following differential inequality
$$
\frac d{dt}\|v(t)\|^p_p\leq p\left(\|\Da
w(t)\|_p+\|w_{xx}(t)\|_p\right)\|v(t)\|_p^{p-1},
$$
which, after integration, leads to
$$
\|v(t)\|_p\leq \|v(t_0)\|_p+\int^t_{t_0}\|\Da
w(\tau)\|_p+\|w_{xx}(\tau)\|_p\,d\tau.
$$
The proof is completed by the result stated in Lemma~\ref{lem 2}.
\end{proof}

\bigskip

Now, we are in a position to prove the main result of this paper

\medskip
\noindent {\it The proof of Theorem \ref{th main}.}
First, we consider the auxiliary solution $\widetilde u=\widetilde
u(x,t)$ of the fractal Burgers
equation \rf{e1.1} with the step-like initial condition \rf{ini-rar}.
In this case,
the measure
$\widetilde m=(u_+-u_-)\delta_0$ is nonnegative, hence
by Theorem \ref{th ux}, $\widetilde u_x\geq 0$ and
by  Lemma \ref{lem 1},
$$
\|u(t)-\widetilde u(t)\|_p\leq Ct^{-(1-1/p)/\alpha}\|u_0-\widetilde
u_0\|_1
$$
for every $p\in [1,\infty]$ and all $t>0$.

Next, we compare $\widetilde u$ with the smooth approximation of the
rarefaction wave
that is  with the solution
$w=w(x,t)$
 of \rf{eq-app}-\rf{ini-app} (note that $\widetilde u_0=w_0$).
By Theorem \ref{th ux} and Lemma
 \ref{le3.6}, we  obtain
$$
\|\widetilde u_x(t)\|_\infty+
\|w_x(t)\|_\infty\leq C t^{-1}.
$$

Moreover, using the following Gagliardo-Nirenberg inequality
\begin{equation}\label{GN}
\|v\|_p\leq C\|v_x\|^a_\infty\|v\|^{1-a}_{p_0},
\end{equation}
 valid
for any $1<p_0<p\leq \infty$
and
$
a=(1/p_0-1/p)/(1+1/p_0),
$
we have
\begin{align*}
\|\widetilde u(t)-w(t)\|_p&
\leq C(\|\widetilde u_x(t)\|_\infty+\|w_x(t)\|_\infty)^a\|\widetilde
u(t)-w(t)\|_{p_0}^{1-a}
\\
&\leq Ct^{-a}\|\widetilde u(t)-w(t)\|^{1-a}_{p_0}.
\end{align*}
Choosing $p_0=(3-\alpha)/(\alpha-1)$ (hence $a=[\alpha-1-(3-\alpha)/p]/2$),
 by Lemma \ref{lem 3}, we conclude $\|\widetilde u(t)-w(t)\|_p\leq
Ct^{-a}\log(2+t)$ for every
$p\in (p_0,\infty]$. Here, we are allowed to use Lemma \ref{lem 3} 
because $\widetilde u_0-w_0\in L^1(\R)\cap L^\infty(\R)\subset 
L^p(\R)$ for every $p\in [1,\infty]$.

Finally, it follows from Lemma \ref{le3.6} that the large time
asymptotics of $w(t)$
is described in $L^p(\R)$ by the rarefaction wave $w^R(x,t)$.

The proof is complete because for $1<\alpha<2$ we have $(1-1/p)/\alpha>(1-1/p)/2$. 
Moreover, since $1<p_0<p$, we have $(1-1/p)/2>(1/p_0-1/p)/(1+1/p_0)$.
\endProof

\section{Additional comments and possible generalizations}
\setcounter{equation}{0}

Our main result is stated and shown in the simplest case of equation
\rf{e1.1}, however, several generalizations are possible.

First of all, the operator $\Da$ can be replaced by the
L\'evy operator ${\cal L}$ which  is a pseudodifferential operator
defined by
the symbol $a=a(\xi)\ge 0$, $\widehat{{\cal L}v}(\xi)=a(\xi)\widehat
v(\xi)$.
Here, the
 function $e^{-ta(\xi)}$ should be positive-definite, so the symbol
$a(\xi)$ can be represented
by the L\'evy--Khintchine formula
 in the Fourier variables
\begin{equation}
a(\xi)=ib\xi+Q(\xi)
+\int_{\R}\bigg(1-e^{-i\eta\xi}
-{i\eta\xi}\,1\!\!\!\;\mbox{I}_{\{|\eta|<1\}}
(\eta)\bigg)\, \Pi(d\eta).\label{L-K}
\end{equation}
Here, $b\in\R$ is fixed, $Q(\xi)=q\xi^2$ with some $q\geq 0$, and 
 $\Pi$ is a Borel measure such that $\Pi(\{0\})=0$ and
$\int_{\R}\min(1,|\eta|^2)\,\Pi(d\eta)<\infty$.

Detailed analysis of conservation laws with the anomalous diffusion
operator
$\cal L$ is contained in papers \cite{bkw,bkw1,bkw2}.
Here, we would like to emphasize that the fundamental nature of the
operator $\cal L$ is clear from the
perspective of probability theory. It represents the most general form
of generator of a~stochastically continuous Markov process with
independent and stationary increments. This fact was our basic
motivation for the development of the theory presented above.

In order to show the convergence toward rarefaction waves of solutions
of conservation laws with the L\'evy operator, we need the
counterparts of estimates
\rf{e3.11}-\rf{e3.12} of the semigroup of linear operators $e^{-t{\cal
L}}$  generated by $-{\cal L}$. They are valid {\it e.g.}
under the assumption
that
the symbol $a$ of ${\cal L}$  has the form
\begin{equation}
a(\xi)=\ell|\xi|^\alpha +k(\xi),
\label{as:2}
\end{equation}
where $\ell>0$, $0<\alpha\leq 2$ and   $k$  is a symbol of another
L\'evy
operator ${\cal K}$ such that
\begin{equation}
\lim_{\xi\to 0}\frac{k(\xi)}{|\xi|^\alpha}=0.\label{as:3}
\end{equation}
The   assumptions \rf{as:2}
and \rf{as:3} are fulfilled, for example,  by {\it
multifractal
diffusion operators\,}
$${\cal L}=-a_0\partial^2_x+\sum_{j=1}^ka_j(-
\partial_x^2)^{\alpha_j/2}$$
with $a_0\ge 0$, $a_j>0$, $1<\alpha_j<2$, and
$\alpha=\min_{1\le j\le k}\alpha_j$.
We refer the reader to \cite{bkw1,bkw2} for the reasoning leading
to the decay estimates of solution of nonlinear problem with operator
$\cal L$ satisfying \rf{as:2}-\rf{as:3}.
That argument can be directly adapted to obtain counterparts of
Theorem \ref{th ux} and Lemma \ref{lem 1} with $\Da$ replaced by
${\cal L}$. Note here that the $L^p$-$L^q$ estimates of the semigroup
$e^{-t{\cal L}}$
are equivalent to a certain Nash inequality, see the papers 
\cite{LS,bkw1}
and the references therein.

Our result also holds true, if we replace the nonlinear term $uu_x$ in
\rf{e1.1}
by $f(u)_x$ with a strictly convex $C^2$-function $f$ (as in the paper
 of
 Il'in
and Oleinik \cite{ilo}) satisfying $f''(u)\geq \kappa$ for some fixed
$\kappa>0$ and all $u\in\R$. Under this assumption, we immediately
generalize Theorem \ref{th ux} and we obtain the decay estimate \rf{ux
dec}.

 In order to show the counterpart of Lemma \ref {lem 1}, we should use
the assumption $f''(u)\geq \kappa$ and to replace equalities \rf{v
pos} by the following (recall that $v=u-\widetilde u$)
$$
\int_\R[f(u)-f(\widetilde u)]_x |v|^{p-2}v\;dx
\geq \kappa \left(1-\frac1p\right)
\int_\R\widetilde u_x|v|^p\;dx\geq 0.
$$

This argument, however, is known and used systematically {\it e.g.} in
\cite[inequality (3.5)]{nin}
(see, also \cite{ilo,man,man1,nak}), hence we skip other details.

\bigskip

{\bf Acknowlegements.}
This paper was partially written  while the first-named
author enjoyed the hospitality and support 
of the Institute of Applied Physics and Computational
Mathematics in Beijing.
G.~Karch was also partially supported 
 by the European Commission Marie Curie Host Fellowship
for the Transfer of Knowledge ``Harmonic Analysis, Nonlinear
Analysis and Probability''  MTKD-CT-2004-013389.
C.~Miao  and X.~Xu  were supported partially  by the NSF of China No.~10571016 and No.~10601009.

\end{document}